\documentclass[12pt]{article}

\usepackage{amsmath}
\usepackage{amsthm}
\usepackage{amssymb}
\usepackage{mathrsfs}
\usepackage{dsfont}
\usepackage{color}
\usepackage{mathabx}

\usepackage[total={6.5in, 9in}]{geometry}
\usepackage[hyperfootnotes=false]{hyperref}

\usepackage[symbol]{footmisc}

\begin{document}

\hypersetup{pdftitle = Argand's "Reflexions'' of 1815 -- An English Translation,
                 pdfauthor = Michael Bertrand,
                 pdfsubject = English translation of Argand's "Reflexions sur la nouvelle théorie des imaginaires'',
                 pdfkeywords = {Jean Robert Argand, History of Complex Numbers, Fundamental Theorem of Algebra}}

\title{Reflexions sur la nouvelle théorie des imaginaires, suivies d'une application à la 
démonstration d'un théorème d’analise\\[0.5in]
	(Reflections on the new theory of imaginaries, followed by an application to the 
demonstration of a theorem of analysis)}
\author{by Jean-Robert Argand\\[0.5in]
       \textit{Annales de Mathématiques} \textbf{5}, No. 7 (January 1, 1815), pp. 197-209\\[0.5in]
	Translated by Michael Bertrand}
\date {October 17, 2022}
\maketitle

\newpage

\LARGE T\normalsize he new theory of imaginaries, which has already been discussed several times in this journal\footnote[1]{Argand's note: See pp. 61, 133, 222 and 364 in the 4th volume.}, has two distinct and independent objects. It tends first of all to give an intelligible meaning to expressions which one was forced to admit in analysis, but which were not believed up to now to be related to any known and evaluable quantity. And secondly, it offers a method of calculation, or, if you will, a \textit{notation} of a particular kind, which employs geometric symbols concurrently with the ordinary algebraic symbols.\\

These two points of view give rise to the two following questions: Does the new theory rigorously demonstrate that $\sqrt{-1}$ is represented by a line perpendicular to the lines taken for +1 and -1? Can the notation of \textit{directed} lines furnish proofs and solutions preferable in simplicity and brevity to those which they seem destined to replace?\\

As for the first point, it will perhaps always be subject to discussion as long as we seek to establish the meaning of $\sqrt{-1}$ \textcolor{red}{[p. 198]} by analogy with the received notions of positive and negative quantities and the proportion between them. We have discussed and will continue to discuss negative quantities, all the more so to rebut objections against the new notions of imaginary numbers.\\

But there will no longer be any difficulty if, as M. Français has done (\textit{Annales}, Vol. IV, p. 62), we establish by definition what we mean by the \textit{ratio} of two lines specified by their \textit{length and position}. Indeed, the relation between two lines of given length and direction is conceived with all necessary geometrical precision. Whether we call this relation a \textit{ratio} or any other name, we can always make it the object of rigorous reasoning and draw from it the geometrical and analytical consequences of which M. Français and I have given a few examples.\\

The only question remaining is therefore whether it is indeed permissible to designate this relation by the words \textit{ratio} or \textit{proportion}, which already have a definite and immutable meaning in analysis. But this is in fact permitted since in the new meaning only \textit{adds to} the old one without \textit{changing} it. One generalizes in such a way that the common understanding is, so to speak, just a particular case of the new approach. It is therefore not a question of looking here for a \textit{demonstration}.\\

It is thus, for example, that the first analyst who said that $a^{-n} = 1/a^n$ gave this equation not as a \textit{theorem} to be proven, but as a \textit{definition} of powers with negative exponents. The only thing he had to show was that by adopting this definition, he was only generalizing the definition of powers with positive exponents, the only ones known until then. It is the same for powers with fractional, irrational or imaginary exponents. It has been said (\textit{Annales}, Vol. IV, p. 231) that Euler demonstrated that ${(\sqrt{-1})}^{\sqrt{-1}} = e^{-1/2 \pi}$. The word \textit{demonstrate} is accurate insofar as one regards this equation as \textcolor{red}{[p. 199]} derived from the equation $e^{x \sqrt{-1}} = \cos{x} + \sqrt{-1} \sin{x}$, from which it easily follows; but it is not comparable to the latter; for to demonstrate that a certain expression has such a value, it is first necessary to have defined this expression; now, do there exist powers with imaginary exponents defined prior to what is called Euler's demonstration? it does not appear so.\\

When Euler sought to reduce the expression $a^{x \sqrt{-1}}$ to evaluable quantities, he naturally considered the previously proven theorem $e^z = 1 + {{z} \over {1}} + {{z^2} \over {1 \cdot 2}} + \cdots$ for all \textit{real} values of $z$. By making $z = x \sqrt{-1}$, he found $e^{x \sqrt{-1}} = 1 + {{x \sqrt{-1}} \over{1}} - {{x^2} \over{2}} - \cdots$; from which he must have concluded, not that $e^{x \sqrt{-1}} = \cos{x} + \sqrt{-1} \sin{x}$, but that, if we defined the expression $e^{x \sqrt{-1}}$ by saying that it represents a quantity equal to $\cos{x} + \sqrt{-1} \sin{x}$, then the powers with real exponents and the powers with imaginary exponents would follow a common pattern. This is therefore only an extension of principles and not the demonstration of a theorem.\\

It is also by an extension of these principles that I was led to look at ${(\sqrt{-1})}^{\sqrt{-1}}$ by representing perpendiculars in the plane by $\pm 1, \pm \sqrt{-1}$. The two approaches are at odds and I am certainly careful not to claim that mine must prevail; I only want to point out that MM. Servois and Français have attacked the problem with considerations which are basically of the same nature as those on which I rely.\\

But if the perpendicular in question cannot be expressed as ${(\sqrt{-1})}^{\sqrt{-1}}$, how will it be expressed? or, to put it better, can we find an expression to represent this perpendicular such that all lines drawn in any direction would be represented, connecting them by a common pattern, so that any line drawn in the plane can be located relative to $\pm1, \pm \sqrt{-1}$? This is a question which seems to excite the curiosity of mathematicians, at least those who accept the new theory.\\

\textcolor{red}{[p. 200]} I return to the first point of discussion by observing that the question whether or not $\sqrt{-1}$ represents a perpendicular to $\pm 1$ brings up the word \textit{ratio}, because everyone understands by this expression a quantity such that $+1 : \sqrt{-1} :: \sqrt{-1} : -1$; that is, the ratios ${{\sqrt{-1}} \over {+1}}, {{-1} \over {\sqrt{-1}}}$ should be equal. Thus the objection made by M. Servois (\textit{Annales}, Vol. IV, p. 228) against the proof of the first theorem of M. Fran\c{c}ais, by saying "that it is not proved that $\pm a \sqrt{-1}$ is a mean of position between $+a$ and $-a$", amounts to saying that the meaning of the word \textit{ratio} does not contain anything relating to position.\\

This is true in the normal sense though even then one could say that the idea of the ratio of two quantities of different signs necessarily brings in those signs. In the new sense, direction joins with magnitude to form a ratio. So it is a simple question of words, which is decided by the precise definition given by Mr. Français, and which is moreover only an extension of the ordinary definition.\\

The second point of discussion is more important. Doubtless there is no truth accessible by the use of the \textit{directed line} notation which cannot also be arrived at by ordinary methods; but will this be achieved more or less easily by one method than by the other? it seems to me that the question deserves to be examined.\\

It is to the influence of new methods and notations in the progressive march of science that the moderns owe their great superiority over the ancients when it comes to mathematical knowledge; thus, when a new idea of this kind presents itself, one should at least examine whether there is any advantage to be drawn from it. Since the publication of the new theory, M. Servois is the only one who has expressed his opinion on this subject, and he does not favor the use of \textit{directed lines} as a notation.\\

The use of analytical formulas seems simpler to him \textcolor{red}{[p. 201]} and more expeditious (\textit{Annales}, Vol. IV, p. 230). I call for a more particular examination of my method. I observe that it is new and that the mental operations which it requires, although very simple, may well require some practice in order to be executed with the dispatch which practice gives to the ordinary operations of algebra. Some of the theorems that I have demonstrated seem to me to be more easily demonstrated than through the purely analytical approach. It may be an author's illusion, and I won't dwell on it; but with some confidence I point to the superiority of directed lines for the proof of the theorem of algebra that ``any polynomial $x^n + ax^{n-1} + \cdots$ can be decomposed into factors of the first or second degree''.\\

I think I should return to this demonstration, both to resolve the objection that M. Servois made to it (\textit{Annales} Vol. IV, p. 231) and to show in more detail how it follows easily from the new approach. The importance and the difficulty of this theorem which has exercised the sagacity of mathematicians of the first order, will excuse in the eyes of the readers, I presume, a few repetitions of what has been said on this same subject.\\

The proofs of this theorem that have been given seem to fall into two classes. Some are based on certain metaphysical principles about functions and the reduction of equations: principles which are doubtless true, but which are not susceptible to rigorous demonstration. They are a species of \textit{axioms}, the truth of which cannot be appreciated unless one already possesses the \textit{spirit} of algebraic calculation; while recognizing the truth of a \textit{theorem} requires understanding the \textit{principles} behind the calculation; that is to say, knowing its definitions and notation. Consequently demonstrations of this kind have been frequently attacked. The journal to which I entrust these reflections, in particular, offers several examples; and the discussions that have taken place on this subject are \textcolor{red}{[p. 202]} an indication that the logic behind them is not entirely without reproach.\\

In other demonstrations, the proposition to be established is attacked head-on by showing that there always exists at least one quantity of the form $a + b \sqrt{-1}$, which taken for $x$, makes the proposed polynomial zero, or that we can solve this polynomial in real factors of the first or second degree. This is the course followed by Lagrange. This great mathematician has shown that the reasoning done before him on this same subject by d'Alembert, Euler, Foncenex, etc., was incomplete (\textit{Résolut. des équat. numériq.} notes IX and X). Some employed expansions in series, others subsidiary equations, but they still had not proven that the coefficients of these equations and of these series were always necessarily real.\\

These mathematicians implicitly admit the principle that "determining an unknown that can be solved in $n$ ways leads to an equation of degree $n$." Lagrange himself regards this as legitimate, although he does not make use of it in the demonstrations cited. Could we not say here as well that this principle, extremely probable no doubt, has not been demonstrated and enters into the class of those sorts of axioms just in questioned. Crucially, although one can acquire belief in this principle by extensive practice in science, this is not the place to employ it, when it is a question of a proposition which is one of the first in the theoretical sphere which presents itself to be demonstrated in analysis. This observation is not intended to raise a quibble on conceptions to which all mathematicians owe the tribute of their esteem, out of place and useless as that would be. It only tends to show the difficulty of dealing with this subject in a satisfactory manner.\\

From these considerations it appears that a demonstration at once direct, simple and rigorous may still merit the attention of mathematicians. I am therefore going to resume here the discussion of page 142 \textcolor{red}{[p. 203]} of volume IV of the \textit{Annales}; and to remove all shadow of doubt, I will free it from the consideration of vanishing quantities.\\

It is worth recalling briefly the first principles of the theory of directed lines. Taking direction $\overline{\mathrm{KA}}$ for positive quantities, the opposite direction $\overline{\mathrm{AK}}$ will as usual indicate negative quantities. Drawing from K the perpendicular BKD, one of the directions $\overline{\mathrm{KB}}, \overline{\mathrm{KD}}$, the first for example, will indicate the imaginaries $+a \sqrt{-1}$, the other the imaginaries $-a \sqrt{-1}$. The line above the letters indicates that the line is drawn from the first point to the second. We suppress the line above the letters when considering only the length of the line.\\

Taking arbitrary points F, G, H, ... P, Q, we have

$$\overline{\mathrm{FG}} + \overline{\mathrm{GH}} + \cdots + \overline{\mathrm{PQ}} = \overline{\mathrm{FQ}}.$$

This is the rule of \textit{addition}.\\

If for four lines we have the equation

$${\mathrm{AB} \over \mathrm{CD}} = {\mathrm{EF} \over \mathrm{GH}},$$

\vspace{10pt}

\noindent and also that the angle between $\overline{\mathrm{AB}}, \overline{\mathrm{CD}}$ is equal to the angle between $\overline{\mathrm{EF}}, \overline{\mathrm{GH}}$, then these lines are said to be in \textit{proportion}. This determines the rule of \textit{multiplication}; for a product is nothing more than a fourth term of proportion where the first is unity.\\

Note that these two rules do not depend on one's opinion of the new theory. If we insist on having $\sqrt{-1}$ as a symbol everywhere in algebra, and which, though sometimes called absurd, has never given results that are; if one wants to say that this symbol is nothing at all, though not zero, that has not caused difficulty. The directed lines will only be the \textit{representations} \textcolor{red}{[p. 204]} of numbers of the form $a + b \sqrt{-1}$.\\

These rules are no less legitimate than the standard operations of algebra, but instead of deducing them \textit{a priori} from partly metaphysical considerations, we will derive the first from a simple construction and the second as an immediate consequence of the formulas $\sin{(a+b)} =\sin{a} \cos{b} + \cdots$, etc.. The benefit is that using these rules results in entirely rigorous demonstrations.\\

Directed lines will therefore signify the numbers $a + b \sqrt{-1}$. Like these numbers, directed lines can be added, subtracted, multiplied, divided, and so on; they mirror the numbers in every respect; in a word, they \textit{represent} them completely. Seen this way, concrete quantities will represent abstract numbers; though the reverse is not true.\\

In what follows, accents indicate the absolute value of the quantities they attach to; so if $a = m + n \sqrt{-1}$ for $m$ and $n$ real numbers, then $a' = \sqrt{m^2 + n^2}$.\\

Let the given polynomial be

$$y(x) = x^n + a x^{n-1} + b x^{n-2} + \cdots + fx + g,$$

\noindent where $n$ is a whole number and $a, b, \cdots f, g$ can be of the form $m + n \sqrt{-1}$. It is a question of proving that one can always find a quantity of this same form which, taken as $x$, results in $y(x)=0$.\\

The polynomial can be evaluated for any value $x$ by the formula above. Taking K as the initial point and P as the final point, the value of the polynomial at $x$ can be be expressed as $\overline{\mathrm{KP}}$ and it must be shown that $x$ can be determined so that the point P coincides with K.\\

Now if, of the infinitly many values which $x$ can take, there were none which result in K and P coinciding, then the line $\overline{\mathrm{KP}}$ \textcolor{red}{[p. 205]} could never become zero; and, of all possible values of $\overline{\mathrm{KP}}$, there would necessarily be one which is smaller than all the others. Let $z$ be the value of $x$ giving this {minimum}; we couldn't have

$$y'(z+i) < y'(z),$$

\noindent for any quantity $i$.\footnote[1]{Translator's note: $i$ is \textit{not} $\sqrt{-1}$, but is a complex number whose value will be determined in what follows.} Expanding gives\\

\noindent $\text{(A)} \hspace{6pt} y(z+i) = y(z) + \{nz^{n-1} + (n-1)\textcolor{red}{a}z^{n-2} + \cdots + f\}i + \left\{{n \over 1} \cdot {{n-1} \over 2}z^2 + \cdots\right\}i^2 + \cdots + \linebreak (nz + a)i^{n-1} + i^n.$\footnote[2]{Translator's note: Argand omits the coefficient $a$ colored red, but this has no bearing on the argument.}\\

Since the coefficients of the different powers of $i$ can be zero, requiring special consideration, it is best to deal with the question in a general way by rewriting (A) as\\

\noindent $\text{(B)} \hspace{6pt} y(z+i) = y(z) + Ri^r + Si^s + \cdots + Vi^v + i^n,$\\

\noindent where none of the coefficients $R, S, ... V$ are zero, and such that the exponents $r, s, \cdots v, n$ increase. Note that if all the coefficients of (A) were zero, then equation (B) would reduce to $y(z+i) = y(z) + i^n$ and then putting $i = \sqrt[n]{-y(z)}$ gives $y(z+i) = 0$, so the theorem would be proven for this case, which can therefore be ignored in what follows. Thus we can assume that the right hand side of equation (B) has at least three terms.\\

Assuming that, we form $y(z+i)$ by taking

$$\overline{\mathrm{KP}} = y(z), \quad \overline{\mathrm{PA}} = Ri^r, \quad \overline{\mathrm{AB}} = Si^s, \; \cdots \; \overline{\mathrm{FG}} = Vi^v, \quad \overline{\mathrm{GH}} = i^n,$$

\noindent so \textcolor{red}{[p. 206]}

$$y'(z) = \mathrm{KP}, \quad R' {(i')}^r = \mathrm{PA}, \quad S' {(i')}^s = \mathrm{AB}, \; \cdots \; V' {(i')}^v = \mathrm{FG}, \quad {(i')}^n = \mathrm{GH},$$

\vspace{10pt}

\noindent because it is clear that in general $p' q' = (pq)'$.\\

$y(z+i)$ will be represented by the broken or straight line $\overline{\mathrm{KPAB} \ldots \mathrm{FGH}}$ or by $\overline{\mathrm{KH}}$ and we have to prove that we can have KH $<$ KP.\\

Now, the quantity $i$ can vary in two ways:\\

\noindent 1) In \textit{direction}: it is obvious that if $i$ has angle $\alpha$, its power $i^r$ will have angle $r \alpha$. Let $\alpha$ be the angle by which $\overline{\mathrm{PA}}=Ri^r$ exceeds $\overline{\mathrm{KP}}=y(z)$. If one gives $i$ to the angle $(\pi - \alpha) / r$, $\overline{\mathrm{PA}}$ will have angle $\pi - \gamma$, that is to say, the direction of $\overline{\mathrm{PA}}$ will be the opposite to that of $\overline{\mathrm{KP}}$; so that the point A will be on the line PK, extended if necessary along its extremity at K.\footnote[1]{Translator's note: The purpose of this paragraph is to show that the direction of $\overline{\mathrm{PA}}=Ri^r$ can always be chosen to be exactly opposite to that of $\overline{\mathrm{KP}}=y(z)$.}\\

\noindent 2) Assuming the direction of $i$ is fixed, we can then vary it by \textit{magnitude}; and first, if PA $>$ KP\footnote[2]{Translator's note: Argand has PA $<$ KP here, but the opposite inequality is clearly meant.}, we can decrease the length of $i$ until PA $<$ KP, so that point A falls between K and P. Then if the magnitude of $i$, thus reduced, is not such that we have

$$R'{i'}^r > S'{i'}^s + \cdots + V'{i'}^v + {i'}^n,$$

\noindent one can, by further diminishing $i$, ensure this inequality because the exponents $s,...v, n$ are all greater than $r$.\\

However, this inequality amounts to

$$\mathrm{PA > AB + \cdots + FG + GH}.$$

The distance AH will therefore be smaller than PA, and consequently if we draw a circle with center A and radius AP, the point H will be inside this circle; it follows from the first elements of geometry that KH $<$ KP, since K is on the extension of the radius PA on the side of the center A.\\

\textcolor{red}{[p. 207]} I invite the reader to draw a figure to follow this demonstration. By applying to it the very simple fundamental principles recalled above we see that with the exception of the expansion (A), which supposes an algebraic calculation, the rest is done on sight with no special effort.\\

It is almost superfluous to pause for an objection that could be made to the preceding, namely, that if one sought to determine the value of $x$ by progressively decreasing $y'(x)$\footnote[1]{Translator's note: $y'(x)$ is not the derivative of $y(x)$, but the absolute value of $y(x)$.} as prescribed, its value might never reach zero because $i'$ might only decrease by smaller and smaller amounts in successive substitutions. In fact it has not been proven that this is impossible, but it just shows that the preceding considerations could not furnish, at least without new developments, a method of approximation; and this in no way invalidates the proof of the theorem.\\

M. Servois' objection is easily resolved. This mathematician seems to say that it is not enough to find values of $x$ which give to the polynomial constantly decreasing values; it is necessary also that the pattern of the decrease brings the polynomial to \textit{zero} -- that zero is not just the \textit{limit} but is actually achieved. But it has been demonstrated that one can find for $y'(x)$ not only ever-decreasing values, but also a value less than that which one would claim to be the smallest of all.\\

If the polynomial cannot be brought to zero, its smallest value will then be other than zero, and in this case the proof retains all its force. M. Servois seems to indicate that he makes a distinction between an infinitely small limit and an absolutely zero limit; if such were the case, one could argue otherwise with considerations quite similar to those that M. Gergonne put forward on an occasion rather analogous to this one; his \textcolor{red}{[p. 208]} answer applies to the present case almost word for word, and, \textit{mutatis mutadis}, it suffices to refer the reader to it (\textit{Annales}, Vol. III, p. 355). Mr. Servois's scruples undoubtedly derive from the consideration of the equation of the hyperbola $y=1/x$. It is indeed certain that, although one can find a value for $y$ in this equation lower than any given limit, $y$ can nevertheless only become zero if one supposes $x$ infinite. But this circumstance does not occur in our demonstration because it is certainly not by an infinite value of $x$ that the polynomial $y'(x)$ becomes zero.\\

Let us return to the subject which gave rise to the issues above; one may ask if it would be possible to translate this approach into the ordinary language of analysis? This seems to me very probable; but perhaps it would be difficult to obtain the result so simply that way. It seems that to do this one would have to bring the expression of the imaginaries closer to the notation of directed lines by writing, for example

$$\sqrt{a^2 + b^2}\left\{{a \over{\sqrt{a^2 + b^2}}} + {b \over{\sqrt{a^2 + b^2}}} \sqrt{-1}\right\} \text{ for } a + b \sqrt{-1}.$$

\noindent $\sqrt{a^2 + b^2}$ could be called the \textit{modulus} of $a+b \sqrt{-1}$ and would represent the \textit{absolute magnitude} of the line $a+b \sqrt{-1}$, while the other coefficient, whose modulus is unity, would represent its direction. One need only prove (1) that \textit{the modulus of the sum of several quantities is not greater than the sum of the moduli of these quantities}, which amounts to saying that the line AF is not greater than the sum of the lines AB, BC, ... EF; and (2) that \textit{the modulus of the product of several quantities is equal to the product of the moduli of these quantities}.\\

I must leave the task of reconciling these methods to more skillful calculators. If one succeeds in this so as to obtain a purely analytical demonstration as simple as that which follows from the new principles, he will have gained something in analysis, thus arriving by an easy route at \textcolor{red}{[p. 209]} a result whose difficulties were worthy of the power of Lagrange himself. If on the contrary one does not succeed in this, the notation of directed lines will then retain an obvious advantage over the ordinary method. In any case, the new theory will have done science a small service.\\

Allow me, in concluding these reflections, to remark on the the note of M. Lacroix inserted in the \textit{Annales} (Vol. IV, p. 367). This learned professor says that the \textit{Philosophical Transactions} of 1806 contain a memoir by M. Buée whose subject is the same as that on which M. Français and I have written. Now it was in this same year, 1806, that I published the \textit{Essai sur une manière de représenter les quantités imaginaires dans les constructions géométriques}, a booklet in which I revealed the principles of the new theory and only an extract of which appeared as a memoir in the 4th volume of the \textit{Annales} (p. 133). We know on the other hand that volumes in academic journals cannot appear until after the year whose date they bear. This should be enough to establish that if, as is quite possible, M. Buée owed only to his own reflections the ideas he developed in his memoir, it still remains certain that I could not have had knowledge of this work when my booklet appeared.

\begin{center}- finis -\end{center}

\end{document}